# INVOLUTIVITY OF DISTRIBUTIONS AT POINTS OF SUPERDENSE TANGENCY WITH RESPECT TO NORMAL CURRENTS

SILVANO DELLADIO

ABSTRACT. Let $\mathcal{D}$ and $T$ be, respectively, a $C^1$ distribution of $k$-planes and a normal $k$-current on $\mathbb{R}^n$. Then $\mathcal{D}$ has to be involutive at almost every superdensity point of the tangency set of $T$ with respect to $\mathcal{D}$.

## 1. INTRODUCTION

Let us consider a distribution $\mathcal{D}$ of $k$-dimensional planes on an open set $\Omega \subset \mathbb{R}^n$ and recall that $\mathcal{D}$ is said to be completely integrable if for each $x \in \Omega$ there exists an integral manifold of $\mathcal{D}$ (that is a $k$- dimensional submanifold $\mathcal{M}$ of $\Omega$ such that the tangent plane to $\mathcal{M}$ at $y$ coincides with $\mathcal{D}(y)$, for each $y \in \mathcal{M}$) through $x$. It is natural to ask under what assumptions on the defining structure, be it a set of differential forms or a set of vector fields, the distribution $\mathcal{D}$ is completely integrable. In the classical context in which it is assumed that $\mathcal{D}$ is of class $C^1$ and the integral manifolds are of class $C^2$, a well-known answer is provided by the following celebrated Frobenius theorem: A distribution is completely integrable if and only if it is involutive at every point of $\Omega$ (cf. Theorem 2.11.9 and Theorem 2.11.11 in [14]). In order to avoid technicalities as much as possible, in this introduction we will not recall the definition of involutive distribution (cf. Section 2.4), but this will not prevent us from giving an idea of the content of this work.

To understand the sense of our main result, we must first point out the following well-known fact, which obviously proves one of the two implications of Frobenius theorem (the easier one): *If $\mathcal{D}$ is of class $C^1$ and $\mathcal{M}$ is an integral manifold of $\mathcal{D}$, then $\mathcal{D}$ is involutive at every point of $\mathcal{M}$.* This property can be generalised through the notion of superdensity. To explain this point, let us consider any $k$-dimensional $C^1$ submanifold $\mathcal{M}$ of $\Omega$ and denote by $\tau(\mathcal{M}, \mathcal{D})$ the tangency set of $\mathcal{M}$ with respect to $\mathcal{D}$, i.e., the set of all points $y \in \mathcal{M}$ such that the tangent plane to $\mathcal{M}$ at $y$ coincides with $\mathcal{D}(y)$. Furthermore, let $\mathcal{H}^k$ be the $k$-dimensional Hausdorff measure in $\mathbb{R}^n$ and let $B_{\mathcal{M}}(x, r)$ be the open metric ball of $\mathcal{M}$ centered at $x \in \mathcal{M}$, of radius $r > 0$ (cf. [4, Section 1.6]). The following property holds (cf. [9, Theorem 1.1]): *Let $x \in \mathcal{M}$ be a $(k+1)$-superdensity point of $\tau(\mathcal{M}, \mathcal{D})$ relative to $\mathcal{M}$, i.e.*

(1.1) $$\mathcal{H}^k(B_\mathcal{M}(x,r) \setminus \tau(\mathcal{M}, \mathcal{D})) = o(r^{k+1}) \qquad (as\ t \to 0+).$$

*Then $x \in \tau(\mathcal{M}, \mathcal{D})$ and $\mathcal{D}$ is involutive at $x$.* This property generalises the "fact" mentioned above. Indeed, if $\mathcal{M}$ is an integral manifold of $\mathcal{D}$ then $\tau(\mathcal{M}, \mathcal{D}) = \mathcal{M}$ and hence (1.1) is trivially satisfied. We observe that this generalisation is equivalent to the following structure result for the tangency set: *If $x \in \mathcal{M}$ and $\mathcal{D}$ is not involutive at $x$, then $x$ is not a $(k+1)$-superdensity point of $\tau(M, \mathcal{D})$ relative to $\mathcal{M}$.* In particular, if $\mathcal{D}$ is nowhere involutive then there are no $(k+1)$-superdensity points of $\tau(\mathcal{M}, \mathcal{D})$ relative to $\mathcal{M}$ (whatever the choice of the $k$-dimensional $C^1$ submanifold $\mathcal{M}$). Despite this, $\tau(\mathcal{M}, \mathcal{D}))$ may be ordinarily dense, i.e., such that $\mathcal{H}^k(B_\mathcal{M}(x,r) \setminus \tau(\mathcal{M}, \mathcal{D})) = o(r^k)$, as $r \to 0+$, for $\mathcal{H}^k$-a.e. $x \in \tau(\mathcal{M}, \mathcal{D})$. In fact, it can be proved that, for every $x \in \Omega$, there exists a $k$-dimensional $C^1$ submanifold $\mathcal{M}$ of $\Omega$ such that $x \in \mathcal{M}$ and $\mathcal{H}^k(\tau(\mathcal{M}, \mathcal{D})) > 0$ (cf. [2]).

One of the main goals of this paper is to prove Corollary 4.1, i.e. a generalisation of [9, Theorem 1.1] in which, instead of $\mathcal{M}$, a normal $k$-current on $\Omega$ is considered. Unfortunately, however, its statement (including the definition of normal $k$-current, cf. Section 2.3 below) is too technical to be used effectively in an introduction such as this, which aims to present the results obtained in a simple and informal way. For the purposes of this presentation, it will be sufficient to simply focus on the application of Corollary 4.1 to integral $k$-currents (which constitute a particularly interesting subfamily of normal $k$-currents). We recall that an integral $k$-current $T$ on $\Omega$ is a linear functional on the space $\mathcal{E}_k$ of smooth and compactly supported differential $k$-forms on $\Omega$, with the following properties:

(i) It is rectifiable with positive integer multiplicity. This means that $T$ is representable by integration as follows

$$\langle T; \omega \rangle = \int_R \langle \eta; \omega \rangle \theta\, d\mathcal{H}^k \qquad \text{(for all } \omega \in \mathcal{E}_k\text{)},$$

where $R$ is a $k$-rectifiable subset of $\Omega$, $\theta$ is a positive integer-valued function in $L^1(\mathcal{H}^k \llcorner R)$ and $\eta$ is a unit simple measurable $k$-vectorfield spanning the approximate tangent $k$-plane to $R$ at $(\mathcal{H}^k \llcorner R)$-a.e. point of $R$.

(ii) The boundary of $T$, that is the $(k-1)$-current $\partial T$ on $\Omega$ defined by

$$\langle \partial T; \omega' \rangle := \langle T; d\omega' \rangle \qquad (\omega' \in \mathcal{E}_{k-1}),$$

is rectifiable with positive integer multiplicity too. Thus there exist a $(k-1)$-rectifiable subset $R'$ of $\Omega$, a positive integer-valued function $\theta' \in L^1(\mathcal{H}^{k-1} \llcorner R')$ and a unit simple measurable $(k-1)$-vectorfield $\eta'$ spanning the approximate tangent $(k-1)$-plane to $R'$ at $(\mathcal{H}^{k-1} \llcorner R')$-a.e. point of $R'$ such that

$$\langle \partial T; \omega' \rangle = \int_{R'} \langle \eta'; \omega' \rangle \theta'\, d\mathcal{H}^k \qquad \text{(for all } \omega' \in \mathcal{E}_{k-1}\text{)}.$$

We now consider a $k$-distribution $\mathcal{D}$ of class $C^1$ on $\Omega$, an integral $k$-current $T$ on $\Omega$ and adopt the notation introduced in (i) and (ii) above. Let us denote by $\Gamma(\eta, \mathcal{D})$ the set of points $x \in R$ at which the approximate tangent $k$-plane to $R$ exists and is equal to $\mathcal{D}(x)$. Moreover let $\Gamma(\eta', \mathcal{D})$ be the set of points $x \in R'$ at which the approximate tangent $(k-1)$-plane to $R'$ exists and is contained in $\mathcal{D}(x)$. Then we have the following result (cf. Corollary 4.2):

**Theorem.** *If $\mathcal{J}$ denotes the set of all $x \in \Omega$ such that*

$$\lim_{r \to 0+} \frac{\int_{B_r(x) \cap (R \setminus \Gamma(\eta, \mathcal{D}))} \theta \, d\mathcal{H}^k}{r^{k+1}} = \lim_{r \to 0+} \frac{\int_{B_r(x) \cap (R' \setminus \Gamma(\eta', \mathcal{D}))} \theta' \, d\mathcal{H}^{k-1}}{r^k} = 0,$$

*then $\mathcal{D}$ is involutive at $(\mathcal{H}^k \llcorner R)$-a.e. $x \in \mathcal{J}$.*

## 2. Basic notation and notions, preliminary results

Throughout this paper $\Omega$ denotes an open subset of $\mathbb{R}^n$ (with $n \geq 2$). The standard basis of $\mathbb{R}^n$ and its dual will be denoted by $e_1, \ldots, e_n$ and $dx_1, \ldots, dx_n$, respectively. If $k$ is any positive integer not exceeding $n$, then $I(n, k)$ is the family of integer multi-indices $\mathbf{i} = (i_1, \ldots, i_k)$ such that $1 \leq i_1 < \cdots < i_k \leq n$. For every $\mathbf{i} = (i_1, \ldots, i_k) \in I(n, k)$, we set

$$e_\mathbf{i} := e_{i_1} \wedge \cdots \wedge e_{i_k}, \quad dx_\mathbf{i} := dx_{i_1} \wedge \cdots \wedge dx_{i_k}.$$

The linear space of $k$-vectors (resp. $k$-covectors) is denoted by $\wedge_k(\mathbb{R}^n)$ (resp. $\wedge^k(\mathbb{R}^n)$). We recall that $\{e_\mathbf{i} \mid \mathbf{i} \in I(n, k)\}$ (resp. $\{dx_\mathbf{i} \mid \mathbf{i} \in I(n, k)\}$) is the standard basis of $\wedge_k(\mathbb{R}^n)$ (resp. $\wedge^k(\mathbb{R}^n)$). Multivectors and multicovectors are in duality. More precisely, the duality between $\wedge_k(\mathbb{R}^n)$ and $\wedge^k(\mathbb{R}^n)$ is defined by

$$\langle \zeta; \alpha \rangle := \sum_{\mathbf{i} \in I(n,k)} \zeta_\mathbf{i} \alpha_\mathbf{i}, \text{ for all } \zeta \in \wedge_k(\mathbb{R}^n) \text{ and } \alpha \in \wedge^k(\mathbb{R}^n),$$

where $\zeta_\mathbf{i}$ (resp. $\alpha_\mathbf{i}$) is the $\mathbf{i}$-th coefficient in the representation of $\zeta$ (resp. $\alpha$) with respect to the standard basis of $\wedge_k(\mathbb{R}^n)$ (resp. $\wedge^k(\mathbb{R}^n)$), that is $\zeta = \sum_{\mathbf{i} \in I(n,k)} \zeta_\mathbf{i} e_\mathbf{i}$ (resp. $\alpha = \sum_{\mathbf{i} \in I(n,k)} \alpha_\mathbf{i} dx_\mathbf{i}$). If $h \leq k$, $\zeta \in \wedge_k(\mathbb{R}^n)$ and $\alpha \in \wedge^h(\mathbb{R}^n)$, then the interior multiplication $\zeta \llcorner \alpha$ is the $(k-h)$-vector defined by

$$\langle \zeta \llcorner \alpha; \beta \rangle = \langle \zeta; \alpha \wedge \beta \rangle, \text{ for all } \beta \in \wedge^{k-h}(\mathbb{R}^n),$$

cf. [10, Section 1.5.1].

The open ball of radius $r$ centered at $x \in \mathbb{R}^n$ is denoted by $B_r(x)$. The Lebesgue measure and the $h$-dimensional Hausdorff measure on $\mathbb{R}^n$ are denoted by $\mathcal{L}^n$ and $\mathcal{H}^h$, respectively. A subset of $\mathbb{R}^n$ is said to be *$h$-rectifiable* if it can be covered, except for an $\mathcal{H}^h$-negligible subset, by countably many $h$-dimensional $C^1$ surfaces. Recall that if $R$ is a $h$-rectifiable subset of $\mathbb{R}^n$, then for $\mathcal{H}^h$-a.e. $x \in R$ there is the approximate tangent $h$-plane to $R$ at $x$ (cf. [13, Theorem 15.19]).

All measures we will consider below (except for $\mathcal{L}^n$ and $\mathcal{H}^h$) will be real-valued and defined on $\mathcal{B}(\Omega)$, that is the $\sigma$-algebra of Borel subsets of $\Omega$. The *restriction* of a measure $\mu$ to $E \in \mathcal{B}(\Omega)$ is defined by

$$(\mu \llcorner E)(B) := \mu(E \cap B), \text{ for all } B \in \mathcal{B}(\Omega).$$

Recall that the *upper $s$-density* and the *lower $s$-density* of $\mu$ at $x \in \Omega$ (with $0 \leq s < +\infty$) are defined by

$$\Theta^s_*(\mu, x) := \liminf_{r \to 0+} \frac{\mu(B_r(x))}{(2r)^s}, \quad \Theta^{*s}(\mu, x) := \limsup_{r \to 0+} \frac{\mu(B_r(x))}{(2r)^s},$$

respectively (cf. [13, Definition 6.8]). Let us also recall the definition of *upper derivative* of another locally finite Borel measure $\lambda$ on $\Omega$ with respect to $\mu$ at $x \in \Omega$:

$$\overline{D}(\lambda, \mu, x) := \limsup_{r \to 0+} \frac{\lambda(B_r(x))}{\mu(B_r(x))},$$

cf. [13, Definition 2.9]. We have the following result (the proof of which is trivial).

**Proposition 2.1.** *Let $\lambda$ and $\mu$ be two locally finite positive Borel measures on $\Omega$. Moreover let $x \in \Omega$ and $s \in [0, +\infty)$ be such that*

$$\Theta^s_*(\mu, x) > 0, \quad \Theta^{*s}(\mu, x) < +\infty.$$

*Then the following inequality holds:*

$$\Theta^s_*(\mu, x) \, \overline{D}(\lambda, \mu, x) \leq \Theta^{*s}(\lambda, x) \leq \Theta^{*s}(\mu, x) \, \overline{D}(\lambda, \mu, x).$$

2.1. **Vectorfields and differential forms.** A map $v : \Omega \to \wedge_k(\mathbb{R}^n)$ is said to be a *$k$-vectorfield*, Analogously, a map $\omega : \Omega \to \wedge^k(\mathbb{R}^n)$ is said to be a *$k$-covectorfield* or (more commonly) a *differential $k$-form*. Obviously, a $k$-vectorfield $v$ (resp. differential $k$-form $\omega$) can be written in terms of the standard basis of $\wedge_k(\mathbb{R}^n)$ (resp. $\wedge^k(\mathbb{R}^n)$), that is,

$$v(x) = \sum_{\mathbf{i} \in I(n,k)} v_\mathbf{i}(x) e_\mathbf{i} \quad \left( \text{resp. } \omega(x) = \sum_{\mathbf{i} \in I(n,k)} \omega_\mathbf{i}(x) dx_\mathbf{i} \right).$$

The regularity of $v$ (resp. $\omega$) is defined by that of its coefficients $v_\mathbf{i}$ (resp. $\omega_\mathbf{i}$). For example, we will say that $v$ (resp. $\omega$) is class $C^1$ if $v_\mathbf{i} \in C^1(\Omega)$ (resp. $\omega_\mathbf{i} \in C^1(\Omega)$) for all $\mathbf{i} \in I(n,k)$.

2.2. **Span of a $k$-vector.** For $v \in \wedge_k(\mathbb{R}^n)$ we define

$$\text{span}(v) := \{v \llcorner \alpha \mid \alpha \in \wedge^{k-1}(\mathbb{R}^n)\}.$$

The span has the following properties (cf. [1, Proposition 5.9]):

(1) if $v = 0$ then $\text{span}(v) = \{0\}$;
(2) if $v \neq 0$ then $\dim \text{span}(v) \geq k$;
(3) if $v_1, \ldots, v_k$ are linearly independent vectors of $\mathbb{R}^n$ and $v = v_1 \wedge \cdots \wedge v_k$, then $\text{span}(v)$ is the $k$-plane generated by $v_1, \ldots, v_k$. In particular $\dim \text{span}(v) = k$;

(4) conversely, if $\dim \mathrm{span}(v) = k$, then $v$ is simple and $v \neq 0$;
(5) $\mathrm{span}(v)$ is the smallest of all linear subspaces $W$ of $\mathbb{R}^n$ such that $v \in \wedge_k(W)$.

We will also need this additional simple property, for which we provide a proof (since we do not have a reference for it).

**Proposition 2.2.** *Let $v \in \wedge_h(\mathbb{R}^n) \setminus \{0\}$ be simple and let $\beta \in \wedge^p(\mathbb{R}^n)$, with $1 \leq p \leq h \leq n$. Assume $v \llcorner \beta = 0$, that is*

(2.1) $$\langle v \llcorner \alpha \, ; \, \beta \rangle = 0, \text{ for all } \alpha \in \wedge^{h-p}(\mathbb{R}^n).$$

*Then $\beta|_{(\mathrm{span}(v))^p} = 0$.*

*Proof.* Consider an orthonormal basis $\varepsilon_1, \ldots, \varepsilon_n$ of $\mathbb{R}^n$ such that $\varepsilon_1, \ldots, \varepsilon_h$ generates $\mathrm{span}(v)$. If $\theta_1, \ldots, \theta_n$ is the dual basis of $\varepsilon_1, \ldots, \varepsilon_n$ and define

$$I_*(n, p) := \{\mathbf{i} = (i_1, \ldots, i_p) \in I(n, p) \,|\, i_p \geq h + 1\},$$

then we have $\langle \varepsilon_\mathbf{i} \, ; \, \beta \rangle$ for all $\mathbf{i} \in I(n, p) \setminus I_*(n, p)$, by (2.1). Then

$$\beta = \sum_{\mathbf{i} \in I_*(n,p)} \langle \varepsilon_\mathbf{i} \, ; \, \beta \rangle \, \theta_\mathbf{i},$$

hence the conclusion follows. $\square$

**Remark 2.1.** *The property established in Proposition 2.2 does not hold, in general, if $v$ is not simple. For example, let $n = 5$, $h = 3$, $p = 2$, $v := e_1 \wedge e_2 \wedge e_3 + e_1 \wedge e_4 \wedge e_5$ and $\beta := dx_2 \wedge dx_4$. Then one can easily prove that the condition (2.1) is verified and that $\mathrm{span}(v) = \mathbb{R}^5$. Since in this case we have $\beta|_{(\mathrm{span}(v))^p} = dx_2 \wedge dx_4$, the above property cannot be valid.*

Consider a Borel map $\tau : \Omega \to \wedge_h(\mathbb{R}^n)$, a Borel differential $l$-form $\omega$ on an open set $U \subset \Omega$ (with $1 \leq l \leq h - 1$) and define the Borel set

$$\mathcal{K}(\tau, \omega) := \{y \in U \,|\, \langle \tau(y) \llcorner \alpha; \omega_y \rangle = 0 \text{ for all } \alpha \in \wedge^{h-l}(\mathbb{R}^n)\}.$$

We observe that in the special case $l = 1$, i.e. if $\omega$ is a Borel differential 1-form, then we have

$$\mathcal{K}(\tau, \omega) = \{y \in U \,|\, \mathrm{span}(\tau(y)) \subset \ker \omega_y\}.$$

2.3. **Currents.** An $h$-*current* on $\Omega$ is a continuous linear functional $T$ on the space $\mathcal{E}_h$ of smooth and compactly supported differential $h$-forms on $\Omega$. The *boundary* of $T$ is an $(h-1)$-current on $\Omega$ denoted with $\partial T$ and defined by $\langle \partial T; \omega \rangle := \langle T; d\omega \rangle$ for every $\omega \in \mathcal{E}_{h-1}$. The *mass* of $T$ is defined as

$$\mathbb{M}(T) := \sup\{\langle T; \omega \rangle \,|\, \omega \in \mathcal{E}_h, \, |\omega(x)| \leq 1 \text{ for every } x \in \Omega\}.$$

Given an $h$-current $T$ on $\Omega$, the following properties are equivalent (by Riesz theorem):

(1) $\mathbb{M}(T) < +\infty$;

(2) There exist a finite positive measure $\mu$ on $\Omega$ and a Borel $h$-vectorfield $\tau$ in $L^1(\mu)$ such that $T = \tau\mu$, i.e.,
$$\langle T; \omega \rangle = \int_\Omega \langle \tau; \omega \rangle \, d\mu \qquad (\omega \in \mathcal{E}_h).$$

Recall from [11, Ch.1, Sect.1.4] that if $\mu$ and $\tau$ are as in (ii) then the total variation of $T = \tau\mu$ equals $|\tau|\mu$, namely
$$|\tau\mu| = |\tau|\mu, \tag{2.2}$$
hence also
$$|\tau\mu|(\Omega) = \int_\Omega |\tau| \, d\mu = \mathbb{M}(T).$$
In particular $|\tau\mu|$ is Radon.

**Remark 2.2.** *Obviously the representation $T = \tau\mu$ is not unique. In particular, we have also $T = \tau\mu \llcorner S_\tau$, with $S_\tau := \{x \in \Omega \,|\, \tau(x) \neq 0\}$. For this reason, it is not restrictive to assume that*
$$\tau(x) \neq 0, \text{ for } \mu\text{-a.e. } x \in \Omega, \tag{2.3}$$
*hence also* $\mathrm{spt}(T) = \mathrm{spt}(\mu)$.

An $h$-current $T$ on $\Omega$ is said to be:

  (i) *normal* if $\mathbb{M}(T)$ and $\mathbb{M}(\partial T)$ are both finite;
  (ii) *rectifiable* if $T = \eta\theta\mathcal{H}^h$ and the following properties hold:
      - $\theta \in L^1(\mathcal{H}^h)$;
      - $R := \{x \in \Omega \,|\, \theta(x) \neq 0\}$ is $k$-rectifiable;
      - $\eta$ is a unit simple $h$-vectorfield such that $\mathrm{span}(\eta(x))$ is the approximate tangent $h$-plane to $R$ at $x$, for $\mathcal{H}^h$-a.e. $x \in R$.
      In this case $T$ is also denoted by $[\![R, \eta, \theta]\!]$.
  (iii) *integral* if $T$ is rectifiable and (with the notation of (ii)):
      - $\theta|_R$ is positive and integer-valued;
      - $\mathbb{M}(\partial T) < +\infty$ (hence $\partial T$).

Recall that if $T$ is integral then also $\partial T$ is integral, cf. [15, Theorem 30.3].

2.4. **Distributions.** A $k$-*distribution* on $\Omega$ (with $1 \leq k \leq n$) is a map $\mathcal{D}$ that to each $x \in \Omega$ associates a $k$-dimensional plane $\mathcal{D}(x) \subset \mathbb{R}^n$. We say that a $k$-distribution $\mathcal{D}$ on $\Omega$ is *of class $C^p$* (with $p \in \mathbb{N}$) if for every $x \in \Omega$ the following property holds: there exist a neighborhood $U \subset \Omega$ of $x$ and a family $\omega^{(1)}, \ldots, \omega^{(n-k)}$ of $C^p$ differential 1-forms on $U$ such that
$$\mathcal{D}(y) = \ker \omega_y^{(1)} \cap \cdots \cap \ker \omega_y^{(n-k)}$$
for all $y \in U$. The forms $\omega^{(1)}, \ldots, \omega^{(n-k)}$ are called *defining forms* (for $\mathcal{D}$) in $U$.

Given an $h$-current with finite mass $T = \tau\mu$ and a $k$-distribution $\mathcal{D}$ on $\Omega$, with $1 \leq h \leq k \leq n$, the tangency set of $T$ with respect to $\mathcal{D}$ is defined as
$$\Gamma(\tau, \mathcal{D}) := \{x \in \Omega \,|\, \mathrm{span}(\tau(x)) \subset \mathcal{D}(x)\}.$$

If $\mathcal{D}$ is a $k$-distribution of class $C^0$ on $\Omega$ and $\omega^{(1)}, \ldots, \omega^{(n-k)}$ are defining forms (for $\mathcal{D}$) in $U \subset \Omega$, then
$$\Gamma(\tau, \mathcal{D}) \cap U = \bigcap_{j=1}^{n-k} \{x \in U \,|\, \mathrm{span}(\tau(x)) \subset \ker \omega_x^{(j)}\},$$

that is

(2.4) $$\Gamma(\tau, \mathcal{D}) \cap U = \bigcap_{j=1}^{n-k} \mathcal{K}(\tau, \omega^{(j)}).$$

**Remark 2.3.** *Let $T = \tau\mu$ be a $k$-current with finite mass in $\Omega$, let $\mathcal{D}$ be a $k$-distribution of class $C^0$ on $\Omega$ and observe that (cf. Section 2.2) the following property holds: If $x \in \Gamma(\tau, \mathcal{D})$ and $\tau(x) \neq 0$, then $\dim \mathrm{span}(\tau(x)) = k$. Hence*
$$\Gamma(\tau, \mathcal{D}) = Z \cup \Gamma_*(\tau, \mathcal{D}),$$

*where*
$$Z := \{x \in \Omega \,|\, \tau(x) = 0\}, \quad \Gamma_*(\tau, \mathcal{D}) := \{x \in \Omega \,|\, \mathrm{span}(\tau(x)) = \mathcal{D}(x)\}.$$

*Observe that*

(2.5) $$|\tau\mu|(Z) = \int_Z |\tau|\, d\mu = 0$$

*by (2.2). Moreover, for all $x \in \Gamma_*(\tau, \mathcal{D})$ the $k$-vector $\tau(x)$ has to be simple (cf. Section 2.2). If we assume the non-restrictive condition (2.3), then (2.5) becomes equivalent to $\mu(Z) = 0$.*

Recall that a $k$-distribution $\mathcal{D}$ of class $C^1$ is said to be *involutive* at $x \in \Omega$ if there exists a family $\omega^{(1)}, \ldots, \omega^{(n-k)}$ of defining forms in a neighborhood of $x$ such that

(2.6) $$(d\omega^{(j)})_x|_{\mathcal{D}(x) \times \mathcal{D}(x)} = 0, \text{ for all } j = 1, \ldots, n-k.$$

One can easily verify that the property (2.6) does not depend on the choice of the family of defining forms. The distribution $\mathcal{D}$ is called *involutive* (in $\Omega$) if it is involutive at every $x \in \Omega$.

Also recall that, if $p \geq 1$ and $\mathcal{D}$ is of class $C^p$, then a non-empty $C^p$ imbedded submanifold $M$ of $\Omega$ such that $\mathcal{T}_x M = \mathcal{D}(x)$ for all $x \in M$ is called a $C^p$ *integral manifold* of $\mathcal{D}$. As a celebrated theorem by Frobenius establishes, the involutivity of $\mathcal{D}$ is a necessary and sufficient condition for the existence of an integral manifold of $\mathcal{D}$ through every point of $\Omega$. This topic is extensively covered in many books of differential geometry, for example in [5, Sect. 3.2], [12, Ch. 19], [14, Sect. 2.11].

## 2.5. Superdensity.
The following definition generalizes the notion of $m$-density point with respect to $\mathcal{L}^n$ (cf. [6, 7, 8]).

**Definition 2.1.** *Let $h \in [0, +\infty)$ and $E \in \mathcal{B}(\Omega)$. Then $x \in \Omega$ is said to be an $h$-superdensity point of $E$ with respect to a Borel measure $\lambda$ if $\lambda(B_r(x) \setminus E) = \lambda(B_r(x)) \, o(r^h)$, as $r \to 0+$. The set of all $h$-superdensity points of $E$ with respect to $\lambda$ is denoted by $E^{\lambda,h}$.*

**Remark 2.4.** *Let $\lambda$ be a Borel measure, $h \in [0, +\infty)$ and $E, F \in \mathcal{B}(\Omega)$. The following facts hold:*

(1) *If $\lambda = \mathcal{L}^n$ then the set of all $h$-superdensity points of $E$ with respect to $\lambda$ coincides with the set of all $(n+h)$-density points of $E$, i.e., $E^{\mathcal{L}^n,h} = E^{(n+h)}$.*
(2) $E^{\lambda,h_2} \subset E^{\lambda,h_1}$, *whenever $0 \leq h_1 \leq h_2 < +\infty$.*
(3) $(E \cap F)^{\lambda,h} = E^{\lambda,h} \cap F^{\lambda,h}$.
(4) $\lambda(E \setminus E^{\lambda,0}) = \lambda(E^{\lambda,0} \setminus E) = 0$ *(cf. [13, Corollary 2.14]).*
(5) *Let $E$ be open. Then $E \subset E^{\lambda,h}$, where the inclusion can be strict (e.g., for $\lambda = \mathcal{L}^n$ and $E = B_r(x) \setminus \{x\}$ one has $E^{\lambda,h} = B_r(x)$). The equality $E = E^{\lambda,h}$ occurs instead in the case that all connected components of $E$ are simply connected.*
(6) *If $x \in \Omega$ and $\lambda(B_r(x)) = 0$ for some $r > 0$, then $x \in E^{\lambda,k}$ for all $k \in [0, +\infty)$.*
(7) $E^{\lambda \llcorner E, k} = \Omega$ *for all $k \in [0, +\infty)$.*

## 3. The main result

Let $h, n$ be integers satisfying $1 \leq h \leq n$ and $\Omega$ be an open subset of $\mathbb{R}^n$. Moreover, let $T$ be a normal $h$-current on $\Omega$. Then $T = \tau \mu$ and $\partial T = \tau' \mu'$, where $\mu, \mu'$ are two finite positive measures on $\Omega$ and
$$\tau : \Omega \to \wedge_h(\mathbb{R}^n), \quad \tau' : \Omega \to \wedge_{h-1}(\mathbb{R}^n)$$
are two Borel maps such that $\tau \in L^1(\mu), \tau' \in L^1(\mu')$, cf. Section 2.3. Recall that

(3.1) $$|\tau \mu| = |\tau|\mu, \quad |\tau' \mu'| = |\tau'|\mu',$$

by (2.2). In particular $|\tau \mu|$ and $|\tau' \mu'|$ are Radon.

**Remark 3.1.** *By [13, Corollary 2.14] there exists $N \subset \Omega$ such that $\mu(N) = 0$ and*

(3.2) $$|\tau(x)| < +\infty, \quad \lim_{r \to 0+} \frac{\int_{B_r(x)} \tau \, d\mu}{\mu(B_r(x))} = \tau(x), \quad \lim_{r \to 0+} \frac{\int_{B_r(x)} |\tau| \, d\mu}{\mu(B_r(x))} = |\tau(x)|$$

*for all $x \in \Omega \setminus N$.*

We also consider a continuous differential $l$-form $\omega$ on an open set $U \subset \Omega$ with $1 \leq l \leq h-1$ and set for simplicity
$$K := \mathcal{K}(\tau, \omega), \quad K' := \mathcal{K}(\tau', \omega).$$

**Remark 3.2.** *We can easily prove that*

$$K^{|\tau\mu|,0} \setminus N \subset K. \tag{3.3}$$

*Indeed, let $x \in K^{|\tau\mu|,0} \setminus N$ and $\alpha \in \wedge^{h-l}(\mathbb{R}^n)$. Then, denoting by $\theta$ the constant differential $(h-l)$-form on $\Omega$ such that $\theta_y = \alpha$ for all $y \in \Omega$, we have*

$$\left| \int_{B_r(x)} \langle \tau; \theta \wedge \omega \rangle \, d\mu \right| = \left| \int_{B_r(x) \setminus K} \langle \tau; \theta \wedge \omega \rangle \, d\mu \right|$$

$$\leq C \int_{B_r(x) \setminus K} |\tau| \, d\mu$$

$$= \left( \int_{B_r(x)} |\tau| \, d\mu \right) o(r^0) \qquad (\text{as } r \to 0+),$$

*by (3.1). Hence*

$$\langle \tau(x) \llcorner \alpha; \omega_x \rangle = \lim_{r \to 0+} \frac{\int_{B_r(x)} \langle \tau; \theta \wedge \omega \rangle \, d\mu}{\mu(B_r(x))} = 0,$$

*by (3.2). Finally (3.3) follows from the arbitrariness of $x \in K^{|\tau\mu|,0} \setminus N$ and $\alpha \in \wedge^{h-l}(\mathbb{R}^n)$.*

**Theorem 3.1.** *Let $T = \tau\mu$, $\partial T = \tau'\mu'$, $\omega$, $K$, $K'$ and $N$ be as above with the additional assumption that $\omega$ is of class $C^1$. Moreover let $x \in K^{|\tau\mu|,1} \setminus N$ and $s \in [0, +\infty)$ be such that*

$$\Theta^s_*(\mu, x) > 0, \quad \Theta^{*s}(\mu, x) < +\infty.$$

*Finally, let $\alpha$ be arbitrarily chosen in $\wedge^{h-l-1}(\mathbb{R}^n)$. Then*

$$|\langle \tau(x) \llcorner \alpha; (d\omega)_x \rangle| \leq C \left( 1 - \frac{\Theta^s_*(\mu, x)}{\Theta^{*s}(\mu, x)} \right) + C \limsup_{r \to 0+} \frac{|\tau'\mu'|(B_r(x) \setminus K')}{r^s}. \tag{3.4}$$

*Proof.* Let $\rho \in (0,1)$ and consider $g \in C^1_c(B_1(0))$ such that $0 \leq g \leq 1$, $g|_{B_\rho(0)} \equiv 1$ and

$$|D_i g| \leq \frac{2}{1-\rho} \qquad (i = 1, \ldots, n).$$

For every real number $r$ such that $0 < r < \text{dist}(x, \mathbb{R}^n \setminus U)$, we define $g_r \in C^1_c(B_r(x))$ as

$$g_r(y) := g\left( \frac{y-x}{r} \right), \quad y \in B_r(x)$$

and observe that (for all $y \in B_r(x)$ and $i = 1, \ldots, n$)

$$|D_i g_r(y)| = \frac{1}{r} \left| D_i g\left( \frac{y-x}{r} \right) \right| \leq \frac{2}{r(1-\rho)}. \tag{3.5}$$

If $\theta$ denotes the constant differential $(h-l-1)$-form on $\Omega$ such that $\theta_y = \alpha$, for all $y \in \Omega$, then

$$d(g_r \, \omega \wedge \theta) = dg_r \wedge \omega \wedge \theta + g_r \, d\omega \wedge \theta,$$

hence

$$\langle \partial T; g_r \, \omega \wedge \theta \rangle = \langle T; dg_r \wedge \omega \wedge \theta + g_r \, d\omega \wedge \theta \rangle,$$

that is
$$\int_\Omega g_r \langle \tau'; \omega \wedge \theta \rangle \, d\mu' = \int_\Omega \langle \tau; dg_r \wedge \omega \wedge \theta \rangle \, d\mu + \int_\Omega g_r \langle \tau; d\omega \wedge \theta \rangle \, d\mu.$$

From now on, for simplicity, we will denote $B_r(x)$ by $B_r$ and $B_{\rho r}(x)$ by $B_{\rho r}$. Recalling the definition of $K$ and $K'$, we obtain

$$\int_{B_r \setminus K'} g_r \langle \tau'; \omega \wedge \theta \rangle \, d\mu' = \int_{B_r \setminus K} \langle \tau; dg_r \wedge \omega \wedge \theta \rangle \, d\mu + \int_{B_r} g_r \langle \tau; d\omega \wedge \theta \rangle \, d\mu$$

and then, by (3.5),

(3.6) $$\left| \int_{B_r} g_r \langle \tau; d\omega \wedge \theta \rangle \, d\mu \right| \leq \frac{C}{r(1-\rho)} \int_{B_r \setminus K} |\tau| \, d\mu + C \int_{B_r \setminus K'} |\tau'| \, d\mu'.$$

On the other hand, we have

(3.7) $$\left| \int_{B_r} g_r \langle \tau; d\omega \wedge \theta \rangle d\mu \right| \geq \left| \int_{B_{\rho r}} g_r \langle \tau; d\omega \wedge \theta \rangle d\mu \right| - \left| \int_{B_r \setminus B_{\rho r}} g_r \langle \tau; d\omega \wedge \theta \rangle d\mu \right|$$
$$= \left| \int_{B_{\rho r}} \langle \tau; d\omega \wedge \theta \rangle \, d\mu \right| - \left| \int_{B_r \setminus B_{\rho r}} g_r \langle \tau; d\omega \wedge \theta \rangle \, d\mu \right|.$$

From (3.6), (3.7) and (3.5) it follows that

$$\left| \int_{B_{\rho r}} \langle \tau; d\omega \wedge \theta \rangle \, d\mu \right| \leq C \int_{B_r \setminus B_{\rho r}} |\tau| \, d\mu + \frac{C}{r(1-\rho)} \int_{B_r \setminus K} |\tau| \, d\mu$$
$$+ C \int_{B_r \setminus K'} |\tau'| \, d\mu',$$

hence (also recalling (3.1))

$$\frac{\mu(B_{\rho r})}{(2\rho r)^s} \cdot \frac{1}{\mu(B_{\rho r})} \left| \int_{B_{\rho r}} \langle \tau; d\omega \wedge \theta \rangle \, d\mu \right| \leq C \left( \frac{\mu(B_r)}{(2\rho r)^s} \cdot \frac{\int_{B_r} |\tau| \, d\mu}{\mu(B_r)} - \frac{\mu(B_{\rho r})}{(2\rho r)^s} \cdot \frac{\int_{B_{\rho r}} |\tau| \, d\mu}{\mu(B_{\rho r})} \right)$$
$$+ \frac{C}{1-\rho} \cdot \frac{\mu(B_r)}{(2\rho r)^s} \cdot \frac{\int_{B_r} |\tau| d\mu}{\mu(B_r)} \cdot \frac{|\tau\mu|(B_r \setminus K)}{r|\tau\mu|(B_r)}$$
$$+ C \frac{|\tau'\mu'|(B_r \setminus K')}{(2\rho r)^s}.$$

Observe that the constant $C$ above is independent from $r$ and $\rho$. Recalling (3.2) and that $x \in K^{|\tau\mu|,1}$, we obtain (letting $r \to 0+$)

$$\Theta^{*s}(\mu, x) \, |\langle \tau(x), (d\omega)_x \wedge \alpha \rangle| = \limsup_{r \to 0+} \frac{\mu(B_{\rho r})}{(2\rho r)^s} \cdot \frac{1}{\mu(B_{\rho r})} \left| \int_{B_{\rho r}} \langle \tau; d\omega \wedge \theta \rangle \, d\mu \right|$$
$$\leq C \left( \rho^{-s} \, \Theta^{*s}(\mu, x) |\tau(x)| - \Theta_*^s(\mu, x) |\tau(x)| \right)$$
$$+ C \, (2\rho)^{-s} \limsup_{r \to 0+} \frac{|\tau'\mu'|(B_r \setminus K')}{r^s}.$$

Thus

$$|\langle \tau(x), (d\omega)_x \wedge \alpha \rangle| \leq C\left(\rho^{-s} - \frac{\Theta^s_*(\mu, x)}{\Theta^{*s}(\mu, x)}\right) + C(2\rho)^{-s} \limsup_{r \to 0+} \frac{|\tau'\mu'|(B_r \setminus K')}{r^s},$$

for all $\rho \in (0, 1)$. The conclusion follows by letting $\rho \to 1-$. □

**Remark 3.3.** *Using Proposition 2.1 with $\lambda = |\tau'\mu'| \llcorner (\Omega \setminus K')$, it is easy to verify that (3.4) is equivalent to the following inequality:*

$$|\langle \tau(x) \llcorner \alpha; (d\omega)_x \rangle| \leq C\left(1 - \frac{\Theta^s_*(\mu, x)}{\Theta^{*s}(\mu, x)}\right) + C \limsup_{r \to 0+} \frac{|\tau'\mu'|(B_r(x) \setminus K')}{\mu(B_r(x))}.$$

## 4. Application to the context of distributions

**Proposition 4.1.** *Let $T = \tau\mu$ be a $k$-current with finite mass in $\Omega$ and let $\mathcal{D}$ be a $k$-distribution of class $C^0$ on $\Omega$. Then there exists $N \subset \Omega$ such that*

$$\mu(N) = 0, \quad \Gamma(\tau, \mathcal{D})^{|\tau\mu|,0} \setminus N \subset \Gamma(\tau, \mathcal{D}).$$

*Proof.* We can find a countable family $B_1, B_2, \ldots$ of open balls of $\mathbb{R}^n$ such that:

(i) $\cup_i B_i = \Omega$;
(ii) for each $i = 1, 2, \ldots$ there exists a family $\omega^{(i,1)}, \ldots, \omega^{(i,n-k)}$ of defining forms for $\mathcal{D}$ in $B_i$ (recall from Section 2.4 that the $\omega^{(i,j)}$ are $C^0$ differential 1-forms on $B_i$).

For all $i = 1, 2, \ldots$ and $j = 1, \ldots, n-k$, define

$$K_{i,j} := \mathcal{K}(\tau, \omega^{(i,j)}) = \{x \in B_i \,|\, \text{span}(\tau(x)) \subset \ker \omega^{(i,j)}_x\}$$

and recall from (2.4), Remark 3.1 and Remark 3.2 that the following facts hold:

(1) $\Gamma(\tau, \mathcal{D}) \cap B_i = \cap_{j=1}^{n-k} K_{i,j}$;
(2) $N_{i,j} \subset \Omega$ has to exist such that

$$\mu(N_{i,j}) = 0, \quad K_{i,j}^{|\tau\mu|,0} \setminus N_{i,j} \subset K_{i,j}.$$

Hence, if define

$$N := \bigcup_{i,j} N_{i,j},$$

we obtain $\mu(N) = 0$ and, for all $i = 1, 2, \ldots$ (by also recalling the properties listed in Remark 2.4),

$$\begin{aligned}
(\Gamma(\tau, \mathcal{D})^{|\tau\mu|,0} \setminus N) \cap B_i &= \Gamma(\tau, \mathcal{D})^{|\tau\mu|,0} \cap B_i \setminus N \\
&= \Gamma(\tau, \mathcal{D})^{|\tau\mu|,0} \cap B_i^{|\tau\mu|,0} \setminus N \\
&= (\Gamma(\tau, \mathcal{D}) \cap B_i)^{|\tau\mu|,0} \setminus N \\
&= \left( \bigcap_{j=1}^{n-k} K_{i,j} \right)^{|\tau\mu|,0} \setminus N \\
&= \bigcap_{j=1}^{n-k} K_{i,j}^{|\tau\mu|,0} \setminus N \\
&\subset \bigcap_{j=1}^{n-k} (K_{i,j}^{|\tau\mu|,0} \setminus N_{i,j}) \\
&\subset \bigcap_{j=1}^{n-k} K_{i,j} \\
&= \Gamma(\tau, \mathcal{D}) \cap B_i.
\end{aligned}$$

The conclusion follows by recalling that the balls $B_i$ cover $\Omega$. $\square$

**Corollary 4.1.** *Let $T$ be a normal $k$-current in $\Omega$, so we have the usual representations $T = \tau\mu$ and $\partial T = \tau'\mu'$ (cf. Section 3). Moreover consider a $k$-distribution $\mathcal{D}$ of class $C^1$ on $\Omega$ and let $\Upsilon$ denote the set of all points $x \in \Omega$ such that:*

(i) $\tau(x) \neq 0$;
(ii) *There exists $s(x) \in [0, +\infty)$ such that $\Theta_*^{s(x)}(\mu, x) = \Theta^{*s(x)}(\mu, x) \in (0, +\infty)$;*
(iii) $x \in \Gamma(\tau, \mathcal{D})^{|\tau\mu|,1}$ *(note that $\Gamma(\tau, \mathcal{D})^{|\tau\mu|,1} = \Gamma_*(\tau, \mathcal{D})^{|\tau\mu|,1}$, by Remark 2.3);*
(iv) $|\tau'\mu'|(B_r(x) \setminus \Gamma(\tau', \mathcal{D})) = o(r^{s(x)})$, *as $r \to 0+$.*

*If $N$ is the $\mu$-null set defined in Proposition 4.1 and $x \in \Upsilon \setminus N$, then the following properties hold:*

(1) *The $k$-vector $\tau(x)$ is simple and $\mathrm{span}(\tau(x)) = \mathcal{D}(x)$;*
(2) *The distribution $\mathcal{D}$ is involutive at $x$.*

*Proof.* (1) We have $\Gamma(\tau, \mathcal{D})^{|\tau\mu|,1} \subset \Gamma(\tau, \mathcal{D})^{|\tau\mu|,0}$, by property (2) in Remark 2.4. Hence $x \in \mathrm{span}(\tau(x)) \subset \mathcal{D}(x)$, by Proposition 4.1. Since $\tau(x) \neq 0$, the conclusion follows from properties (2) and (4) in Section 2.2.

(2) Let $\{B_i\}$, $\{\omega^{(i,j)}\}$, $\{K_{i,j}\}$ and $\{N_{i,j}\}$ be the families defined in the proof of Proposition 4.1 (here we can obviously assume that the $\omega^{(i,j)}$ are of class $C^1$). Without loss of generality we can suppose that $x \in B_1$. Then, by recalling assumption (iii), properties (3) and (5) in Remark 2.4 and (2.4), we obtain

$$
\begin{aligned}
(4.1) \quad x \in \Upsilon \cap B_1 \setminus N &\subset \Gamma(\tau, \mathcal{D})^{|\tau\mu|,1} \cap B_1 \setminus N \\
&= \Gamma(\tau, \mathcal{D})^{|\tau\mu|,1} \cap B_1^{|\tau\mu|,1} \setminus N \\
&= (\Gamma(\tau, \mathcal{D}) \cap B_1)^{|\tau\mu|,1} \setminus N \\
&= \left( \bigcap_{j=1}^{n-k} K_{1,j} \right)^{|\tau\mu|,1} \setminus N \\
&= \bigcap_{j=1}^{n-k} K_{1,j}^{|\tau\mu|,1} \setminus N \\
&\subset \bigcap_{j=1}^{n-k} \left( K_{1,j}^{|\tau\mu|,1} \setminus N_{1,j} \right).
\end{aligned}
$$

Moreover (by (2.4))

$$\Gamma(\tau', \mathcal{D}) \cap B_1 = \bigcap_{j=1}^{n-k} K'_{1,j},$$

where

$$K'_{1,j} := \mathcal{K}(\tau', \omega^{(1,j)}) \qquad (j = 1, \ldots, n-k).$$

Hence

$$B_r(x) \setminus K'_{1,j} \subset B_r(x) \setminus \Gamma(\tau', \mathcal{D}) \qquad (j = 1, \ldots, n-k),$$

provided $r$ is small enough. Recalling also (iv), we obtain

$$(4.2) \qquad |\tau'\mu'|(B_r(x) \setminus K'_{1,j}) = o(r^{s(x)}) \qquad (j = 1, \ldots, n-k)$$

as $r \to 0+$. Now (ii), (4.1), (4.2) and Theorem 3.1 yield

$$\langle \tau(x) \llcorner \alpha; (d\omega^{(1,j)})_x \rangle = 0 \qquad (j = 1, \ldots, n-k)$$

for all $\alpha \in \wedge^{k-2}(\mathbb{R}^n)$. From Proposition 2.2 we obtain

$$(d\omega^{(1,j)})_x|_{\text{span}(\tau(x)) \times \text{span}(\tau(x))} = 0 \qquad (j = 1, \ldots, n-k).$$

Now the conclusion follows from statement (1). □

**Remark 4.1.** *Let $\mathcal{M}$ be a closed $k$-dimensional $C^1$ submanifold of $\Omega$ with $C^1$ boundary such that $\mathcal{H}^k(\mathcal{M})$ and $\mathcal{H}^{k-1}(\partial\mathcal{M})$ are finite. Let $\tau_\mathcal{M}$ and $\tau'_\mathcal{M}$ be, respectively, a continuous unit simple $k$-vectorfield orienting $\mathcal{M}$ and a continuous unit simple $(k-1)$-vectorfield orienting $\partial\mathcal{M}$ such that the Stoke's identity*

$$\int_\mathcal{M} \langle \tau_\mathcal{M}; d\omega \rangle \, d\mathcal{H}^k = \int_{\partial\mathcal{M}} \langle \tau'_\mathcal{M}; \omega \rangle \, d\mathcal{H}^{k-1}$$

holds for all $C^1$ differential $(k-1)$-forms with compact support in $\Omega$. Then we consider the maps $\tau : \Omega \to \wedge_k(\mathbb{R}^n)$ and $\tau' : \Omega \to \wedge_{k-1}(\mathbb{R}^n)$ extending $\tau_{\mathcal{M}}$ and $\tau'_{\mathcal{M}}$, respectively, such that $\tau|_{\Omega\setminus\mathcal{M}} \equiv 0$ and $\tau'|_{\Omega\setminus\partial\mathcal{M}} \equiv 0$. We observe that:

(1) $T := \tau \mathcal{H}^k \llcorner \mathcal{M}$ is a normal $k$-current, with $\partial T = \tau' \mathcal{H}^{k-1} \llcorner \partial \mathcal{M}$.
(2) The equations (3.2) hold for all $x \in \mathcal{M}$, hence we can assume that the set $N$ introduced in Section 3 coincides with $\Omega \setminus \mathcal{M}$.

Now set for simplicity
$$\{\mathrm{span}(\tau_{\mathcal{M}}) = \mathcal{D}\} := \{y \in \mathcal{M} \,|\, \mathrm{span}(\tau_{\mathcal{M}}(y)) = \mathcal{D}(y)\} = \Gamma(\tau, \mathcal{D}) \bigcap \mathcal{M},$$
$$\{\mathrm{span}(\tau'_{\mathcal{M}}) \subset \mathcal{D}\} := \{y \in \partial\mathcal{M} \,|\, \mathrm{span}(\tau'_{\mathcal{M}}(y)) \subset \mathcal{D}(y)\} = \Gamma(\tau', \mathcal{D}) \bigcap \partial\mathcal{M}$$
and let us consider
$$x \in \{\mathrm{span}(\tau_{\mathcal{M}}) = \mathcal{D}\}^{\mathcal{H}^k \llcorner \mathcal{M}, 1} \bigcap \{\mathrm{span}(\tau'_{\mathcal{M}}) \subset \mathcal{D}\}^{\mathcal{H}^{k-1} \llcorner \partial\mathcal{M}, 1} \bigcap \mathcal{M}.$$
We observe that
$$\{\mathrm{span}(\tau_{\mathcal{M}}) = \mathcal{D}\}^{\mathcal{H}^k \llcorner \mathcal{M}, 1} = \Gamma(\tau, \mathcal{D})^{\mathcal{H}^k \llcorner \mathcal{M}, 1} \bigcap \mathcal{M}^{\mathcal{H}^k \llcorner \mathcal{M}, 1} = \Gamma(\tau, \mathcal{D})^{\mathcal{H}^k \llcorner \mathcal{M}, 1}$$
by (3) and (7) in Remark 2.4. Analogously,
$$\{\mathrm{span}(\tau'_{\mathcal{M}}) \subset \mathcal{D}\}^{\mathcal{H}^{k-1} \llcorner \partial\mathcal{M}, 1} = \Gamma(\tau', \mathcal{D})^{\mathcal{H}^{k-1} \llcorner \partial\mathcal{M}, 1} \bigcap \partial\mathcal{M}^{\mathcal{H}^{k-1} \llcorner \partial\mathcal{M}, 1}$$
$$= \Gamma(\tau', \mathcal{D})^{\mathcal{H}^{k-1} \llcorner \partial\mathcal{M}, 1}.$$
Hence we easily obtain
$$x \in \Upsilon \bigcap \mathcal{M} = \Upsilon \setminus N.$$
Now, by applying Corollary 4.1, we conclude that $\mathrm{span}(\tau_{\mathcal{M}}(x)) = \mathcal{D}(x)$ and $\mathcal{D}$ is involutive at $x$.

The following corollary generalises the property established in Remark 4.1 for the smooth case.

**Corollary 4.2.** *Let $\mathcal{D}$ and $T$ be, respectively, a $k$-distribution of class $C^1$ in $\Omega$ and an integral $k$-current in $\Omega$. Moreover, if $T = [\![R, \eta, \theta]\!]$ and $\partial T = [\![R', \eta', \theta']\!]$, let $\mathcal{J}$ be the set of all $x \in \Omega$ such that*

$$(4.3) \quad \lim_{r \to 0+} \frac{\int_{B_r(x) \cap (R\setminus\Gamma(\eta, \mathcal{D}))} \theta \, d\mathcal{H}^k}{r^{k+1}} = \lim_{r \to 0+} \frac{\int_{B_r(x) \cap (R'\setminus\Gamma(\eta', \mathcal{D}))} \theta' \, d\mathcal{H}^{k-1}}{r^k} = 0.$$

*Then $\mathcal{D}$ is involutive at $(\mathcal{H}^k \llcorner R)$-a.e. $x \in \mathcal{J}$.*

*Proof.* Let us define
$$\mu := \mathcal{H}^k \llcorner R, \quad \tau := \theta\eta, \quad \mu' := \mathcal{H}^{k-1} \llcorner R', \quad \tau' := \theta'\eta'$$

so that
$$T = \tau\mu, \quad \partial T = \tau'\mu'.$$

From [13, Corollary 2.14] and [13, Theorem 17.6] it follows that the following equalities hold at $\mu$-a.e. $x \in \Omega$:

$$\text{(4.4)} \qquad \lim_{r \to 0+} \frac{\int_{B_r(x)} \theta \, d\mu}{\mu(B_r(x))} = \theta(x), \quad \lim_{r \to 0+} \frac{\mu(B_r(x))}{(2r)^k} = 1,$$

hence also

$$\text{(4.5)} \qquad \lim_{r \to 0+} \frac{\int_{B_r(x)} \theta \, d\mu}{r^k} = 2^k \theta(x).$$

We shall prove the thesis by applying Corollary 4.1. To this end, it will suffice to prove that conditions (i-iv) of Corollary 4.1 are verified at $\mu$-a.e. $x \in \mathcal{J}$ (that is $\mu(\mathcal{J} \setminus \Upsilon) = 0$), which we do below:

- Assumption (i) is verified at $\mu$-a.e. $x \in \Omega$, since $T$ is rectifiable (cf. Section 2.3).
- Assumption (ii) is verified at $\mu$-a.e. $x \in \Omega$ (with $s(x) = k$), by the second equality of (4.4).
- Since $\Gamma_*(\tau, \mathcal{D}) = \Gamma(\eta, \mathcal{D})$, we have

$$\frac{|\tau\mu|(B_r(x) \setminus \Gamma_*(\tau, \mathcal{D}))}{|\tau\mu|(B_r(x))} = \frac{\int_{B_r(x) \cap (R \setminus \Gamma(\eta, \mathcal{D}))} \theta \, d\mathcal{H}^k}{\int_{B_r(x)} \theta \, d\mu}$$
$$= \frac{\int_{B_r(x) \cap (R \setminus \Gamma(\eta, \mathcal{D}))} \theta \, d\mathcal{H}^k}{r^{k+1}} \cdot \left(\frac{\int_{B_r(x)} \theta \, d\mu}{r^k}\right)^{-1} r.$$

  Hence, recalling also (4.3) and (4.5), we find that assumption (iii) is verified at $\mu$-a.e. $x \in \mathcal{J}$.
- If define $Z' := \{x \in \Omega \,|\, \tau'(x) = 0\}$, then we have

$$\Gamma(\tau', \mathcal{D}) = Z' \bigcup \Gamma(\eta', \mathcal{D}), \quad |\tau'\mu'|(Z') = 0.$$

Thus
$$|\tau'\mu'|(B_r(x) \setminus \Gamma(\tau', \mathcal{D})) = |\tau'\mu'|(B_r(x) \setminus \Gamma(\eta', \mathcal{D}))$$
$$= \int_{B_r(x) \cap (R' \setminus \Gamma(\eta', \mathcal{D}))} \theta' \, d\mathcal{H}^{k-1}.$$

From this equality and (4.3) it follows that assumption (iv) is verified at every $x \in \mathcal{J}$.

□

University of Trento, Department of Mathematics, via Sommarive 14, 38123 Trento, Italy

*Email address*: `silvano.delladio@unitn.it`